\documentclass[11pt]{article}
\usepackage[applemac]{inputenc}
\usepackage{amsmath,amsthm,amsfonts,amssymb}
\usepackage[a4paper,body={150mm,220mm},centering]{geometry}
\usepackage{fancyhdr}
\usepackage[small]{titlesec}
\newcommand{\dotafter}[1]{#1.}
\titleformat{\section}[hang]
{\normalfont\Large\bf\sffamily}{\thesection.}{.5em}{\dotafter}[]
\titleformat{\subsection}[runin]
{\normalfont\large\bf\sffamily}{\thesubsection.}{.4em}{}[.]
\titlespacing*{\subsection}{0pt}{3ex plus 1ex minus .2ex}{.9em}
\newtheorem{thm}{Theorem}
\newtheorem{lemme}[thm]{Lemma}
\newtheorem{prop}[thm]{Proposition}

\theoremstyle{definition}

\theoremstyle{remark}

\newcommand{\rset}{\mathbb{R}}
\newcommand{\qset}{\mathbb{Q}}

\newcommand{\nset}{\mathbb{N}}
\newcommand{\nl}{\nolimits}

\newcommand{\ind}{\mathbf{1}}
\newcommand{\fl}{\longrightarrow}
\newcommand{\e}{\mathbb{E}}
\newcommand{\p}{\mathbb{P}}
\newcommand{\lp}{\mathrm{L}}
\newcommand{\m}{\mathcal}
\DeclareMathOperator{\sgn}{sgn}
\begin{document}
\title{\bf\sffamily BSDE with quadratic growth and unbounded terminal value}
\author{Philippe Briand \and Ying Hu}
\date{\textit{IRMAR, Université Rennes 1, 35 042 RENNES Cedex, FRANCE} \\[1ex]
\texttt{\normalsize philippe.briand@univ-rennes1.fr \quad ying.hu@univ-rennes1.fr}} 
\maketitle
\thispagestyle{fancy}
\begin{quote}
\textbf{\sffamily Abstract.}\quad In this paper, we study the existence of solution to BSDE with quadratic growth and unbounded terminal value. We apply a localization procedure together with a priori bounds. As a byproduct, we apply the same method to extend a result on BSDEs with integrable terminal condition.
\end{quote}
\section{Introduction}
In this paper we are concerned with real valued backward stochastic differential equations -- BSDEs for short in the remaining --
$$
Y_t=\xi + \int_t^T f(s,Y_s,Z_s)\,ds-\int_t^T Z_s\cdot dB_s,\qquad 0\leq t\leq T
$$
where $(B_t)_{t\geq 0}$ is a standard brownian motion. Such equations have been extensively studied since the first paper of E.~Pardoux and S.~Peng~\cite{PP90}. The full list of contributions is too long to give and we will only quote results in our framework.

Our setting is mainly the following : the generator, namely the function $f$, is of quadratic growth in the variable $z$ and the terminal condition, the random variable $\xi$, will not be bounded. BSDEs with quadratic growth have been first  studied by Magdalena~Kobylanski in her PhD (see \cite{Kob97, Kob00}) and then by Jean-Pierre Lepeltier and Jaime San Martin in~\cite{LsM98}. We should point out that BSDEs with quadratic growth in the variable $z$ have found applications in control and finance, see, e.g., Bismut~\cite{Bis78}, El Karoui, Rouge \cite{ER00}, Hu, Imkeller, Muller \cite{HIM05}, \ldots 

All the results on BSDEs with quadratic growth require that the terminal condition $\xi$ is a bounded random variable. The boundedness of the terminal condition appears, from the point of view of the applications,  to be restrictive and, moreover, from a theoretical point of view, is not necessary to obtain a solution. Indeed, let us consider the following well known equation :
$$
Y_t = \xi + \frac{1}{2} \int_t^T |Z_s|^2\,ds -\int_t^T Z_s\,dB_s,\qquad 0\leq t\leq T~;
$$
the change of variables $P_t = e^{Y_t}$, $Q_t = e^{Y_t} Z_t$, leads to the equation
$$
P_t = e^{\xi} -\int_t^T Q_s\,dB_s
$$
which has a solution as soon as $e^{\xi}$ is integrable. 

On this simple example we see that the existence of exponential moments of the terminal condition is sufficient to construct a solution to our BSDE. Our paper will be focused on the theoretical study of these BSDEs but with unbounded terminal value with only exponential moments. 

To fill the gap between boundedness and existence of exponential moments, we  will use an approach based upon a localization procedure together with a priori bounds. Let us quickly explain how it works on a simple example. Let $f:\rset\times\rset^d\fl\rset$ be a continuous function and $\xi$ be a nonnegative terminal condition such that
$$
|f(y,z) | \leq \frac{1}{2} \, |z|^2,\qquad \e\left[ e^\xi \right] < \infty,
$$
and let us try to construct a solution to the BSDE
$$
Y_t = \xi +\int_t^T f(Y_s,Z_s)\,ds -\int_t^T Z_s\cdot dB_s,\quad 0\leq t\leq T.
$$
As mentioned before, BSDEs with quadratic growth in the variable $z$ can be solved when the terminal solution is bounded. That is why we introduce $(Y^n,Z^n)$ as the minimal solution to the BSDE
$$
Y^n_t = \xi\wedge n +\int_t^T f(Y^n_s,Z^n_s)\,ds -\int_t^T Z^n_s\cdot dB_s,
$$
and of course we want to pass to the limit when $n\to\infty$ in this equation.

The process $Y^n$ is known to be bounded but the estimate depends on $\| \xi\wedge n \|_\infty$ and thus is far from being useful when $\xi$ is not bounded. The first step of our approach consists in finding an estimation for $Y^n$ independent of $n$. In this example, we can use the explicit formula mentioned before to show that
$$
0 \leq -\ln \e\left( e^{-(\xi\wedge n)} \:|\: \m F_t \right) \leq Y^n_t \leq \ln \e\left( e^{\xi\wedge n} \:|\: \m F_t \right) \leq \ln \e\left( e^\xi \:|\: \m F_t \right).
$$

With these inequalities in hands, we introduce the stopping time 
$$
\tau_k = \inf\left\{t\in[0,T] : \ln \e\left( e^\xi \:|\: \m F_t \right) \geq k \right\}\wedge T
$$
and instead of working on the time interval $[0,T]$ we will restrict ourselves to $[0,\tau_k]$ by considering the BSDE
$$
Y^n_{t\wedge \tau_k} = Y^n_{\tau_k} +\int_{t\wedge\tau_k}^{T\wedge\tau_k} f\left( Y^n_s,Z^n_s\right) ds - \int_{t\wedge\tau_k}^{T\wedge\tau_k} Z^n_s\cdot dB_s,\qquad 0\leq t\leq T.
$$

By construction, we have $\sup_n \sup_t \left\| Y^n_{t\wedge \tau_k} \right\|_\infty \leq k$. This last property together with the fact that the sequence $(Y^n)_{n\geq 1}$ is nondecreasing allows us, with the help of a result of Kobylanski, to pass to the limit when $n\to\infty$, $k$ being fixed and then to send $k$ to infinity to get a solution.

\bigskip

The rest of the paper is organized as follows. Next section is devoted to the notations we use during this text. In Section~\ref{quad}, we claim our main result that we prove in Section~\ref{quadproof}. Section~\ref{comp} is devoted to some additional results on BSDEs with quadratic growth in $z$. Finally, in the last section, we apply the same approach to study BSDEs with terminal value in $\lp^1$.

\section{Notations}
For the remaining of the paper, let us fix a nonnegative real number $T>0$.

First of all, $B=\{B_t\}_{t\geq 0}$ is a standard brownian
motion with values in $\rset^d$ defined on some complete probability space
$(\Omega,\mathcal{F},\p)$. $\{\mathcal{F}_t\}_{t\geq 0}$ is the
augmented natural filtration of $B$ which satisfies the usual
conditions. In this paper, we will always use this filtration. $\m P$ denotes the sigma-field of predictable subsets of
$[0,T]\times\Omega$.

As mentioned in the introduction, a BSDE is an equation of the following type
\begin{equation}
\label{maineq}
Y_t=\xi + \int_t^T f(s,Y_s,Z_s)\,ds-\int_t^T Z_s\cdot dB_s,\qquad 0\leq t\leq T.
\end{equation}
$f$ is called the generator and $\xi$ the terminal condition.

Let us recall that a generator is  a random function
$f:[0,T]\times\Omega\times\rset\times\rset^{d}\fl\rset$ which is measurable with respect to $\m P\times\mathcal{B}(\rset)\times\mathcal{B}(\rset^d)$ and a terminal condition is simply a real $\mathcal{F}_T$--measurable random variable.

By a solution to the BSDE~\eqref{maineq} we mean a pair $(Y,Z)=\{ (Y_t,Z_t)\}_{t\in[0,T]}$ of predictable processes  with values in $\rset\times\rset^{d}$ 
such that $\p$--a.s., $t\longmapsto Y_t$ is continuous, $t\longmapsto Z_t$ belongs to  $\lp^2(0,T)$, $t\longmapsto f(t,Y_t,Z_t)$ belongs to $\lp^1(0,T)$
and $\p$--a.s.  
$$
Y_t = \xi +\int_t^T f(s,Y_s,Z_s)\, ds -\int_t^T Z_s\cdot dB_s, \qquad
0\leq t\leq T.
$$

We will use the notation BSDE($\xi,f$) to say that we consider the BSDE whose generator is $f$ and whose terminal condition is $\xi$; $\left(Y^f(\xi),Z^f(\xi) \right)$ means a solution to the BSDE($\xi,f$).

$\left(Y^f(\xi),Z^f(\xi) \right)$ is said to be minimal if $\p$-a.s., for each $t\in[0,T]$, $ Y^f_t(\xi) \leq Y^g_t(\zeta)$ whenever $\p$--a.s. $\xi \leq \zeta$ and $f(t,y,z)\leq g(t,y,z)$ for all $(t,y,z)$. $\left(Y^f(\xi),Z^f(\xi) \right)$ is said to be minimal in some space $\m B$ if it belongs to this space and the previous property holds true as soon as $\left(Y^g(\zeta),Z^g(\zeta) \right)\in\m B$.

\medskip

For any real $p>0$, $\mathcal{S}^p$ denotes the set of real-valued, adapted and c\`adl\`ag processes
$\{Y_t\}_{t\in[0,T]}$ such that 
$$
\left\| Y \right\|_{\m S^p} :=\e\left[\sup\nl_{t\in[0,T]} |Y_t|^p\right]^{1\wedge 1/p} < +\infty.
$$
If $p\geq 1$, $\|\cdot\|_{\mathcal{S}^p}$ is a norm on
$\mathcal{S}^p$ and if 
$p\in(0,1)$, $(X,X')\longmapsto \big\|X-X'\big\|_{\mathcal{S}^p}$ defines a distance on $\mathcal{S}^p$. Under this metric, $\mathcal{S}^p$ is
complete. 

$\mathrm{M}^p$ denotes the set of (equivalent classes of)
predictable processes $\{Z_t\}_{t\in[0,T]}$ with values
in $\rset^d$ such that
$$
\left\| Z \right\|_{\mathrm{M}^p} : = \e\left[\Big(\int_0^T |Z_s|^2\,
ds\Big)^{p/2}\right]^{1 \wedge 1/p} < +\infty. 
$$ 
For $p\geq 1$, $\mathrm{M}^p(\rset^n)$ is a Banach space endowed with
this norm  and for $p\in(0,1)$, $\mathrm{M}^p$ is a complete metric
space with the resulting distance. 

We set $\m S = \cup_{p>1} \m S^p$, $\mathrm{M} = \cup_{p>1} \mathrm{M}^p$ and denote by $\m S^\infty$ the set of predictable bounded processes. Finally, let us recall that a continuous process $\{Y_t\}_{t\in[0,T]}$ belongs to the class (D) if the family $\{ Y_\tau : \tau \mbox{ stopping time bounded by } T \}$ is uniformly integrable.

\section{Quadratic BSDEs}
\label{quad}
In this section, we consider BSDE($\xi,f$) when the generator $f$ has a linear growth in $y$ and a quadratic growth in $z$. We denote \eqref{h1} the assumption: there exist $\alpha \geq 0$, $\beta \geq 0$ and $\gamma>0$  such that $\p$--a.s.
\begin{equation}
\label{h1}
\tag{H1}
\begin{split}
& \forall t\in[0,T], \qquad (y,z)\longmapsto f(t,y,z) \mbox{ is continuous},  \\
&\forall (t,y,z)\in[0,T]\times\rset\times\rset^d,\qquad  \left| f(t,y,z) \right| \leq \alpha + \beta |y| + \frac{\gamma}{2} |z|^2 . 
\end{split}
\end{equation}

Concerning the terminal condition $\xi$, we will assume that \begin{equation}
\label{h2}
\tag{H2}
\e\left[ e^{\gamma e^{\beta T}\, |\xi| }\right] < +\infty.
\end{equation}
We will use also a stronger assumption on the integrability of $\xi$ namely
\begin{equation}
\label{h3}
\tag{H3}
\exists \lambda > \gamma e^{\beta T},\qquad \e\left[ e^{\lambda\, |\xi| }\right] < +\infty.
\end{equation}

It is clear that we can assume without loss of generality that $\alpha \geq \beta / \gamma$.

\medskip

As we explained in the introduction, our method relies heavily on a priori estimate. To obtain such estimations, we will use the change of variable $P_t= e^{\gamma Y_t}$, $Q_t = \gamma e^{\gamma Y_t} Z_t$; if $(Y,Z)$ is a solution to the BSDE($\xi,f$), $(P,Q)$ solves the BSDE
$$
P_t = e^{\gamma \xi} + \int_t^T F(s,P_s,Q_s)\, ds - \int_t^T Q_s\cdot dB_s,\quad 0\leq t\leq T,
$$
with the function $F$ defined by
\begin{equation}
\label{capf}
F(s,p,q) = \ind_{p>0} \left(\gamma p\, f\left(s,\frac{\ln p}{\gamma},\frac{q}{\gamma p}\right) - \frac{1}{2} \, \frac{|q|^2}{p}\right).
\end{equation} 

In view of the growth of the generator $f$, we have $F(s,p,q) \leq \ind_{p>0}\, p (\alpha\gamma + \beta |\ln p| )$. For notational convenience, we denote by $H$ the function
$$
\forall p \in\rset,\qquad H (p) = p \left( \alpha\gamma  + \beta \ln p \right) \ind_{[1,+\infty)} (p) + \gamma\alpha \ind_{(-\infty,1)}(p).
$$
It is straightforward to check that, since $\alpha \geq \beta/\gamma$, $H$ is convex and locally Lipschitz continuous and that, for any real $p>0$, $p \left(\alpha\gamma + \beta |\ln p |  \right) \leq H(p)$.
Thus we deduce the inequality
\begin{equation}
\label{triv}
\forall s\in[0,T],\;\forall p\in\rset,\;\forall q\in \rset^d,\qquad F(s,p,q) \leq H(p).
\end{equation}

To get an upper bound for $Y_t$, the idea is to compare more or less $P_t$ with $\phi_t(\xi)$ where, for any real $z$, $\{ \phi_t(z) \}_{0\leq t\leq T}$ stands for the solution to the differential equation 
\begin{equation}
\label{eqode}
\phi_t = e^{\gamma z} + \int_t^T H(\phi_s)\, ds,\qquad 0\leq t\leq T.  
\end{equation}
Using the convexity of $H$, we will able to prove that
$$
P_t \leq \e\left( \phi_t(\xi) \:|\:\m F_t \right),\qquad Y_t \leq \frac{1}{\gamma} \, \ln \e\left( \phi_t(\xi) \:|\:\m F_t\right).
$$

Before proving this result rigorously, let us recall that the differential equation~\eqref{eqode} can be solved easily. Indeed, we have, for any $z\geq 0$, 
$$
\phi_t(z) = \exp\left(\gamma\alpha\frac{e^{\beta(T-t)}-1}{\beta}\right) \exp\left(z\gamma e^{\beta(T-t)}\right),\quad\mbox{ if }\beta >0,
$$
and $\phi_t(z) = e^{\gamma\alpha (T-t)} e^{\gamma z}$ if $\beta = 0$. 

Let us consider the case where $z<0$. If $e^{\gamma z} + T\gamma\alpha \leq 1$ then the solution is 
$$
\phi_t = e^{\gamma z} +  \gamma\alpha (T-t)
$$ 
and otherwise there exists $ 0<S<T$ such that $e^{\gamma z} +  \gamma\alpha (T-S) = 1$ and 
$$
\phi_t = \left[ e^{\gamma z} +  \gamma\alpha (T-t)\right] \ind_{(S,T]}(t) + \exp\left(\gamma\alpha\frac{e^{\beta(S-t)}-1}{\beta}\right)\, \ind_{[0,S]}(t).
$$
It is plain to check that $t\mapsto \phi_t(z)$ is decreasing and that $z\mapsto \phi_t(z)$ is increasing and continuous.

\begin{lemme}
\label{estbdd}
Let the assumption~\eqref{h1} hold and let $\xi$ be a bounded $\m F_T$--measurable random variable. 

If $(Y,Z)$ is a solution to the BSDE($\xi,f$) in $\m S^\infty\times\mathrm{M}^2$ then
$$
-\frac{1}{\gamma} \ln \e\left(\phi_t(-\xi) \:|\: \m F_t \right)\leq Y_t \leq \frac{1}{\gamma} \ln \e\left(\phi_t(\xi)\:|\: \m F_t \right).
$$
\end{lemme}

\begin{proof}
Let us set $\Phi_t = \e\left( \phi_t(\xi) \:|\: \m F_t\right)$. We have 
$$
\Phi_t = \e\left( e^{\gamma \xi} + \int_t^T H(\phi_s(\xi)) \, ds \:\Big|\: \m F_t \right) = \e\left( e^{\gamma \xi} + \int_t^T \e\left( H(\phi_s(\xi)) \:|\: \m F_s\right)\, ds \:\Big|\: \m F_t \right).
$$
Thus writing the bounded brownian martingale
$$
\e\left( e^{\gamma \xi} + \int_0^T \e\left( H(\phi_s(\xi)) \:|\: \m F_s\right)\, ds \:\Big|\: \m F_t \right) = \e\left[ e^{\gamma \xi} + \int_0^T \e\left( H(\phi_s(\xi)) \:|\: \m F_s\right)\, ds  \right] + \int_0^t \Psi_s\cdot dB_s
$$
$(\Phi,\Psi)$ solves the BSDE
$$
\Phi_t = e^{\gamma \xi} + \int_t^T \e\left( H(\phi_s(\xi)) \:|\: \m F_s\right)\, ds  - \int_t^T \Psi_s\cdot dB_s.
$$
On the other hand, if $(Y,Z)\in \m S^\infty\times\mathrm{M}^2$ is a solution of \eqref{maineq}, setting as before $P_t = e^{\gamma Y_t}$, $Q_t = \gamma e^{\gamma Y_t} Z_t$, we have
$$
P_t = e^{\gamma \xi} + \int_t^T F(s,P_s,Q_s)\, ds - \int_t^T Q_s\cdot dB_s,
$$
with $F$ defined by \eqref{capf}.

It follows that
$$
\Phi_t-P_t = \int_t^T \left(H(\Phi_s)-H(P_s)\right) ds + \int_t^T R_s \,ds - \int_t^T (\Psi_s-Q_s)\cdot dB_s
$$
where, in view of the inequality \eqref{triv} and since $H$ is convex,
$$
R_s = \e\left( H(\phi_s(\xi)) \:|\: \m F_s\right) - H\left(\e\left( \phi_s(\xi) \:|\: \m F_s\right)\right) + H(P_s) - F(s,P_s,Q_s) 
$$
is a nonnegative process.

$H$ is only locally Lipschitz but since $\Phi$ and $P$ are bounded  we can apply the comparison theorem to get $P_t \leq \Phi_t$ and $Y_t \leq \frac{1}{\gamma} \ln \Phi_t$.

Finally, since the function $-f(t,-y,-z)$ still satisfies the assumption \eqref{h1}, we get also the inequality $-Y_t \leq \frac{1}{\gamma} \ln \e\left( \phi_t(-\xi)\:|\: \m F_t\right)$. 
\end{proof}

We are now in position to prove that under the assumptions described before the BSDE~\eqref{maineq} has at least a solution.

\begin{thm}
\label{mainthm}
Let the assumptions \eqref{h1} and \eqref{h2} hold. Then the BSDE~\eqref{maineq} has at least a solution $(Y,Z)$ such that :
\begin{equation}
\label{mainest}
 -\frac{1}{\gamma} \ln \e\left(\phi_t(-\xi) \:|\: \m F_t \right)\leq Y_t \leq \frac{1}{\gamma} \ln \e\left(\phi_t(\xi)\:|\: \m F_t \right).
\end{equation}

If moreover, \eqref{h3} holds, then $Z$ belongs to $\mathrm{M}^2$.
\end{thm}

\begin{proof}[Proof of the last part of Theorem~\ref{mainthm}]
If $(Y,Z)$ is a solution to the BSDE~\eqref{maineq} such that the inequalities~\eqref{mainest} hold, then
$$
|Y_t| \leq \frac{1}{\gamma} \ln \e\left(\phi_0(|\xi|)\:|\: \m F_t \right)
$$
and, under the assumption~\eqref{h3}, we deduce that, for some $p>1$,
$$
\e\left[ \sup\nolimits_{t\in[0,T]} e^{p\gamma |Y_t|} \right] < +\infty.
$$  

For $n\geq 1$, let $\tau_n$ be the following stopping time
$$
\tau_n = \inf\left\{ t\geq 0 : \int_0^t e^{2\gamma |Y_s|} |Z_s|^2 \, ds \geq n \right\} \wedge T,
$$
and let us consider the function from $\rset_+$ into itself defined by
$$
u(x) = \frac{1}{\gamma^2} \left( e^{\gamma x} - 1 - \gamma x \right).
$$
$x\longmapsto u(|x|)$ is $\m C^2$ and we have from Itô's formula, with the notation $\sgn(x) = - \ind_{x\leq 0} + \ind_{x>0}$,
\begin{eqnarray*}
u(|Y_0|) & = & u(|Y_{t\wedge\tau_n}|) + \int_0^{t\wedge\tau_n} \left( u'(|Y_s|) \sgn(Y_s) f(s,Y_s,Z_s) - \frac{1}{2} u''(|Y_s|) |Z_s|^2 \right) ds  \\
& & - \int_0^{t\wedge\tau_n} u'(|Y_s|) \sgn(Y_s) Z_s\cdot dB_s.
\end{eqnarray*}
It follows from \eqref{h1} since $u'(x) \geq 0$ for $x\geq 0$ that
\begin{eqnarray*}
u(|Y_0|) & \leq & u(|Y_{t\wedge\tau_n}|) + \int_0^{t\wedge\tau_n} u'(|Y_s|) \left(\alpha  + \beta |Y_s|\right) ds - \int_0^{t\wedge\tau_n} u'(|Y_s|) \sgn(Y_s) Z_s\cdot dB_s \\
& & - \frac{1}{2} \int_0^{t\wedge\tau_n} \left(u''(|Y_s|) - \gamma\, u'(|Y_s|) \right) |Z_s|^2 ds .
\end{eqnarray*}

Moreover, we have $(u'' - \gamma u')(x) = 1$ for $x\geq 0$ and, taking expectation of the previous inequality, we get 
$$
\frac{1}{2}\, \e\left[\int_0^{T\wedge\tau_n} |Z_s|^2\,ds\right] \leq \e\left[\frac{1}{\gamma^2}\,\sup_{t\in[0,T]} e^{\gamma |Y_t|} + \frac{1}{\gamma}\, \int_0^T  e^{\gamma |Y_s|} \left( \alpha  + \beta |Y_s| \right) ds\right]
$$

Fatou's lemma together with the fact that $e^{\gamma |Y_t|}\in \m S^p$ gives the result.
\end{proof}

\section{Proof of Theorem~\ref{mainthm}}
\label{quadproof}
Let us first construct a solution to the BSDE~\eqref{maineq} in the case where $\xi$ is nonnegative.

For each $n\in\nset^*$, we set $\xi^n = \xi\wedge n$. Then it is known from \cite[Theorem 2.3]{Kob00} that the BSDE
$$
Y^{n}_t = \xi^{n} + \int_t^T f\left(s,Y^{n}_s , Z^{n}_s \right) ds - \int_t^T Z^{n}_s\cdot dB_s, \qquad 0\leq t \leq T
$$
has a minimal solution $(Y^n,Z^n)$ in $\m S^\infty \times \mathrm{M}^2$.  Lemma~\ref{estbdd} implies the inequalities
$$
-\frac{1}{\gamma} \ln \e\left(\phi_t\left(-\xi^{n}\right) \:|\: \m F_t \right)\leq Y^{n}_t \leq \frac{1}{\gamma} \ln \e\left(\phi_t\left(\xi^{n}\right)\:|\: \m F_t \right) .
$$ 

Since we consider only minimal solutions, we have, 
$$
\forall t\in[0,T],\qquad  Y^{n}_t \leq Y^{n+1}_t.
$$
We define $Y = \sup_{n\geq 1} Y^{n}$.

Since $0\leq \phi_t(\xi^n) \leq \phi_0(|\xi|)$ and $0\leq \phi_t(-\xi^n) \leq \phi_0(|\xi|)$, we deduce from the dominated convergence theorem, noting that the random variable $\phi_0(|\xi|)$ is integrable by $\eqref{h2}$, that
$$
-\frac{1}{\gamma} \ln \e\left(\phi_t\left(-\xi\right) \:|\: \m F_t \right)\leq Y_t \leq \frac{1}{\gamma} \ln \e\left(\phi_t\left(\xi\right)\:|\: \m F_t \right).
$$

In particular, we have $\lim_{t\to T} Y_t = \xi =Y_T$. Indeed, for each $S<T$,
$$
\limsup_{t\to T} Y_t \leq \limsup_{t\to T} \frac{1}{\gamma} \ln \e\left(\phi_t\left(\xi\right)\:|\: \m F_t \right) \leq \lim_{t\to T}\frac{1}{\gamma} \ln \e\left(\phi_S\left(\xi\right)\:|\: \m F_t \right) = \frac{1}{\gamma} \ln \phi_S(\xi),
$$
and $\lim_{S\to T} \frac{1}{\gamma} \ln \phi_S(\xi) = \xi$. We can do the same for $\liminf$.

Let us introduce the following stopping time :
$$
\tau_k =\inf\left\{ t\in[0,T] :  \frac{1}{\gamma} \ln \e\left(\phi_0\left(|\xi|\right) \:|\: \m F_t \right) \geq k \right\}\wedge T.
$$

Then $(Y^n_k,Z^n_k):=(Y^{n}_{t\wedge \tau_k},Z^{n}_t\ind_{t\leq \tau_k})$ satisfies the following BSDE 
$$
Y^n_k = \xi^n_k + \int_t^T \ind_{s\leq \tau_k} f\left(s,Y^{n}_k(s),Z^{n}_k(s)\right) ds - \int_t^T Z^{n}_k(s)\cdot dB_s,
$$
where of course $\xi^n_k = Y^n_k(T) = Y^n_{\tau_k}$. 

We are going to pass to the limit when $n$ tends to $+\infty$ for $k$ fixed in this last equation. 

The key point is that $Y^{n}_k$ is increasing in $n$ and remains bounded by $k$.  At this stage, let us mention a mere generalization of  Proposition 2.4 in \cite{Kob00}.

\begin{lemme}[\cite{Kob00}]
\label{kobm}
Let  $(\xi_n)_{n\geq 1}$ be a sequence of $\m F_T$--measurable bounded random variables and $(f_n)_{n\geq 1}$ be a sequence  of generators which are continuous with respect to $(y,z)$. 

We assume that $(\xi_n)_{n\geq 1}$ converges $\p$--a.s. to $\xi$, that $(f_n)_{n\geq 1}$ converges locally uniformly in $(y,z)$ to the generator $f$, and also that 
\begin{enumerate}
\item $\sup_{n\geq 1} \| \xi_n\|_\infty < +\infty$ ;
\item $\sup_{n\geq 1} \left| f_n(t,y,z)\right|$ satisfies the inequality in \eqref{h1}.
\end{enumerate}

If for each $n\geq 1$, the BSDE($\xi_n,f_n$) has a solution in $\m S^\infty\times\mathrm{M}^2$, such that
 $\left( Y^{f_n}(\xi_n) \right)_{n\geq 1}$ is nondecreasing (respectively nonincreasing),
then $\p$--a.s. $\left(Y^{f_n}_t(\xi_n)\right)_{n\geq 1}$ converges uniformly on $[0,T]$ to $Y_t=\sup_{n\geq 1} Y^{f_n}_t(\xi_n)$ (respectively $Y_t=\inf_{n\geq 1} Y^{f_n}_t(\xi_n)$),  $\left(Z^{f_n}(\xi_n)\right)_{n\geq 1}$ converges to some $Z$ in $\mathrm{M}^2$ and $(Y,Z)$ is a solution to BSDE($\xi,f$) in $\m S^\infty\times\mathrm{M}^2$. 
\end{lemme}

\begin{proof}
It follows from Lemma~\ref{estbdd} that there exists $r>0$ such that, $\p$--a.s.
$$
\forall n\geq 1,\quad \forall t \in[0,T], \qquad \left|Y^{f_n}_t(\xi_n) \right| \leq r.
$$

Let us consider the continuous function $\rho(x) = x r / \max (r, |x|)$. Since $\rho(x) = x$ for $|x|\leq r$, $\left(Y^{f_n}(\xi_n),Z^{f_n}(\xi_n) \right)$ solves the BSDE($\xi_n,g_n$) where $g_n(t,y,z) = f_n(t,\rho(y),z)$. Obviously, we have, for each $n\geq 1$,
$$
\left| g_n(t,y,z) \right | \leq \alpha + \beta\, r + \frac{\gamma}{2} |z|^2,
$$
and thus we can apply the result of Kobylanski.
\end{proof}

Setting $Y_k(t) = \sup_n Y^{n}_k(t)$, it follows from the previous lemma that there exists a process $Z_k\in \mathrm{M}^2$ such that $\lim_n Z^{n}_k = Z_k$ in $\mathrm{M}^2$ and $(Y_k, Z_k)$ solves the BSDE
\begin{equation}
\label{edsrk}
Y_k(t) = \xi_k  + \int_t^T \ind_{s\leq \tau_k} f\left(s,Y_{k}(s),Z_{k}(s)\right) ds - \int_t^T Z_{k}(s) \cdot dB_s,
\end{equation}
where $\xi_k = \sup_n Y^{n}_{\tau_k}$.

But $\tau_k\leq \tau_{k+1}$, and thus we get, coming back to the definition of $Y_k$, $Z_k$ and $Y$,
$$
Y_{t\wedge\tau_k} = Y_{k+1}(t\wedge\tau_k)=Y_k(t), \qquad 
Z_{k+1}(t)\,\ind_{t\leq \tau_k}=Z_{k}(t).
$$

As $\tau_k\rightarrow T$ and the $Y_k$'s are continuous processes we deduce  in particular that $Y$ is continuous on $[0,T)$. On the other hand, as mentioned before $\lim_{t\rightarrow T} Y_t=\xi$ and $Y_T$ is equal to $\xi$ by construction. Thus $Y$ is a continuous process on the closed interval $[0,T]$.

Then we define $Z$ on $(0,T)$ by setting :
$$
Z_t=Z_k(t), \quad \mbox{if } t\in (0,\tau_k).
$$

From \eqref{edsrk}, $(Y,Z)$ satisfies:
\begin{equation}
\label{edsr2}
Y_{t\wedge\tau_k}=Y_{\tau_k}+\int_{t\wedge\tau_k}^{\tau_k} f(s,{Y}_s,{Z}_s)ds-\int_{t\wedge\tau_k}^{\tau_k}Z_s\cdot dB_s.
\end{equation}

Finally, we have
\begin{eqnarray*}
\p(\int_0^T|Z_s|^2\, ds=\infty)
&= & \p\left(\int_0^T|Z_s|^2\,ds=\infty,\tau_k=T\right)+\p\left(\int_0^T|Z_s|^2\,ds=\infty,\tau_k<T\right)\\
&\leq & \p\left(\int_0^{\tau_k}|Z_k(s)|^2\,ds=\infty\right)+\p(\tau_k<T),
\end{eqnarray*}
and we deduce that, $\p$--a.s.
$$
\int_0^T |Z_s|^2\,ds<\infty.
$$

By sending $k$ to infinity in \eqref{edsr2}, we deduce that $(Y,Z)$ is a solution of \eqref{maineq}.

\medskip

Let us explain quickly how to extend this construction to the general case. Let us fix $n\in\nset^*$ and $p\in\nset^*$ and set $\xi^{n,p} = \xi^+\wedge n - \xi^- \wedge p$. Let us consider, $(Y^{n,p},Z^{n,p})$ the minimal bounded solution to the BSDE
$$
Y^{n,p}_t = \xi^{n,p} + \int_t^T f\left(s,Y^{n,p}_s , Z^{n,p}_s \right) ds - \int_t^T Z^{n,p}_s\cdot dB_s, \qquad 0\leq t \leq T
$$
which satisfies 
$$
-\frac{1}{\gamma} \ln \e\left(\phi_t\left(-\xi^{n,p}\right) \:|\: \m F_t \right)\leq Y^{n,p}_t \leq \frac{1}{\gamma} \ln \e\left(\phi_t\left(\xi^{n,p}\right)\:|\: \m F_t \right) .
$$

We have,
$$
\forall t\in[0,T],\qquad Y^{n,p+1}_t \leq Y^{n,p}_t \leq Y^{n+1,p}_t,
$$
and we define $Y^p = \sup_{n\geq 1} Y^{n,p}$ so that $Y^{p+1}_t \leq Y^p_t$ and $Y_t = \inf_{p\geq 1} Y^p_t$.

By the dominated convergence theorem, we have
$$
-\frac{1}{\gamma} \ln \e\left(\phi_t\left(-\xi\right) \:|\: \m F_t \right)\leq Y_t \leq \frac{1}{\gamma} \ln \e\left(\phi_t\left(\xi\right)\:|\: \m F_t \right),
$$
and in particular, we have $\lim_{t\to T} Y_t = \xi =Y_T$. 

$(Y^{n,p}_{t\wedge \tau_k},Z^{n,p}_t\ind_{t\leq \tau_k})$ solves the BSDE 
$$
Y^{n,p}_{t\wedge \tau_k} = Y^{n,p}_{\tau_k} + \int_t^T \ind_{s\leq \tau_k} f\left(s,Y^{n,p}_s,Z^{n,p}_s\right) ds - \int_t^T Z^{n,p}_s \ind_{s\leq \tau_k} \cdot dB_s.
$$
But, once again $Y^{n,p}_{t\wedge\tau_k}$ is increasing in $n$ and decreasing in $p$ and remains bounded by $k$. Arguing as before, setting $Y_k(t) = \inf_p\sup_n Y^{n,p}_{t\wedge \tau_k}$, there exists a process $Z_k$ such that $\lim_p\lim_n Z^{n,p}(s) \ind_{s\leq \tau_k} = Z_k(s)$ and $(Y_k, Z_k)$ still solves the BSDE~\eqref{edsrk}. The rest of the proof is unchanged.

\section{Additional results on quadratic BSDEs}
\label{comp}

\subsection{Minimal solution}
In this section, we give some complements on BSDEs with quadratic growth in $z$.
\begin{prop}
Let \eqref{h1} holds and assume moreover that there exists an integer $r\geq 0$ such that $\p$--a.s.
$$
f(t,y,z) \geq -r\left(1+|y|+|z|\right).
$$
Let us assume also that \eqref{h3} holds for $\xi^+$ and that, for some $p>1$, $\xi^-\in\lp^p$.

Then BSDE($\xi,f$) has a minimal solution in $\m S$.
\end{prop}

\begin{proof}
For each $n\geq r$, let us consider the function 
$$
f_n(t,y,z) = \inf\left\{ f(t,p,q) + n |p-y| + n |q-z| : (p,q) \in \qset^{1+d} \right\}.
$$
Then $f_n$ is well defined and it is globally Lipschitz continuous with constant $n$. Moreover $(f_n)_{n\geq r}$ is increasing and converges pointwise to $f$. Dini's theorem implies that the convergence is also uniform on compact sets. We have also, for all $n\geq r$,
$$
 -r(1+|y| + |z|) \leq f_n(t,y,z) \leq f(t,y,z)
$$

Let $(Y^n,Z^n)$ be the unique solution in $\m S^p\times \mathrm{M^p}$ to BSDE($\xi,f_n$). It follows from the classical comparison theorem that 
$$
Y^r_t \leq Y^n_t \leq Y^{n+1}_t.
$$

Let us prove that $Y^n_t \leq \frac{1}{\gamma} \ln \e\left( \phi_t(\xi) \:|\: \m F_t\right)$. To do this let us recall that, since $f_n$ is Lipschitz, $Y^n_t = \lim_{m\to +\infty} Y^{f_n}_t(\xi_m)$ where $\xi_m = \xi\,\ind_{|\xi| \leq m}$.  Moreover $\e\left( \phi_t(\xi_m) \:|\: \m F_t\right) \fl \e\left( \phi_t(\xi) \:|\: \m F_t\right)$ a.s. since $\sup_{m\geq 1} |\phi_t(\xi_m)| \leq \phi_0(\xi^+)$ which is integrable. Thus we have only to prove that $Y^{f_n}_t(\xi_m) \leq \frac{1}{\gamma} \ln \e\left( \phi_t(\xi_m) \:|\: \m F_t\right)$.  We keep the notations of the beginning of Section~\ref{quadproof}. $(\Phi, \Psi)$ is solution to the BSDE
$$
\Phi_t = e^{\gamma \xi_m} + \int_t^T H(\Phi_s) \, ds + \int_t^T \Gamma_s\,ds - \int_t^T \Psi_s\cdot dB_s,
$$
where $\Gamma_s =  \e\left( H(\phi_s(\xi) ) \:|\: \m F_s\right) - H(\Phi_s) $ is a nonnegative process since $H$ is convex.

It follows by setting $U_t = \frac{1}{\gamma} \ln \Phi_t$, $V_t = \frac{\Psi_t}{\gamma \Phi_t}$ that $(U,V)$ solves the BSDE
$$
U_t = \xi_m + \int_t^T g(s,U_s,V_s)\, ds  -\int_t^T V_s\cdot dB_s 
$$
where we have set $ g(s,u,v) = (\alpha+\beta u)\ind_{u\geq 0} + \alpha e^{\gamma |u|} \ind_{u<0} +\frac{\gamma}{2} |v|^2 + C_s$ with $C_s = \frac{1}{\gamma} e^{-\gamma U_s} \Gamma_s$. Since the process $C$ is still nonnegative, we have the inequalities 
$$
f_n(t,u,v)\leq f(t,u,v) \leq g(t,u,v)
$$ 
taking into account the fact that $\alpha\gamma \geq \beta$.

Since $f_n$ is Lipschitz continuous and $\left( Y^{f_n}(\xi_m) - U\right)^+$ belongs to $\m S$, we can apply the extended comparison theorem (see Proposition~\ref{wcomp}) to get, for each $m\geq 1$, $Y^{f_n}_t(\xi_m) \leq U_t$ and thus the inequality we want to obtain.

We set $Y = \sup_{n \geq r} Y^n$ and, for $k\geq 1$, 
$$
\tau_k = \inf\left\{ t\in[0,T] :  \max\left(\frac{1}{\gamma} \ln \e\left( \phi_t(\xi) \:|\: \m F_t\right),-Y^r_t\right) \geq k \right\}\wedge T.
$$
Arguing as in the proof of Theorem~\ref{mainthm}, we construct a process $Z$ such that $(Y,Z)$ solves BSDE($\xi,f$). 

Let us show that this solution is minimal in $\m S$. Let $(Y',Z')$ be a solution to the BSDE($\xi',f'$) where $\xi\leq \xi'$ and $f\leq f'$. It is enough to check that $Y^n \leq Y'$ to prove that $Y \leq Y'$. But this is a direct consequence of Proposition~\ref{wcomp}.
\end{proof}

To be complete, let us claim and prove the extended comparison theorem that we used in the proof of the previous result.

\begin{prop}
\label{wcomp}
Let $(Y,Z)$ be a solution to BSDE($\xi,f$) and $(Y',Z')$ be a solution to BSDE($\xi',f'$). We assume that $\xi\leq \xi'$ and that $f$ satisfies, for some constants $\mu$ and $\lambda$, $\p$--a.s.
\begin{eqnarray*}
&&(y-y')\cdot \left( f(t,y,z)-f(t,y',z) \right) \leq \mu |y-y'|^2; \\
&&\left| f(t,y,z) - f(t,y,z') \right| \leq \lambda |z-z'|; 
\end{eqnarray*}

If $\left(Y-Y'\right)^+$ belongs to $\m S$, then $\p$--a.s. $Y_t \leq Y'_t$.
\end{prop}

\begin{proof}
Let us fix $n\in\nset^*$ and denote $\tau_n$ the stopping time
$$
\tau_n = \inf\left\{ t\in[0,T] : \int_0^t \left(|Z_s|^2 + \left|Z'_s\right|^2\right) ds \geq n \right\} \wedge T.
$$

Tanaka's formula leads to the equation, setting $U_t = Y_t-Y'_t$, $V_t=Z_t-Z'_t$, 
\begin{equation}
\label{tan}
\begin{split}
e^{\mu(t\wedge\tau_n)} U^+_{t\wedge\tau_n} \leq & e^{\mu\tau_n} U^+_{\tau_n}    - \int_{t\wedge\tau_n}^{\tau_n} e^{\mu s} \ind_{U_s>0} V_s\cdot dB_s  \\
&  + \int_{t\wedge\tau_n}^{\tau_n} e^{\mu s} \left\{ \ind_{U_s >0}\left( f(s,Y_s,Z_s) - f'\left(s,Y'_s,Z'_s\right) \right) - \mu U^+_s \right\} ds .
\end{split}
\end{equation} 

First of all, we write 
$$
f(s,Y_s,Z_s) - f'\left(s,Y'_s,Z'_s\right) = f(s,Y_s,Z_s) - f\left(s,Y'_s,Z_s\right) + f\left(s,Y'_s, Z_s\right) -f'\left(s,Y'_s,Z'_s\right)
$$
and we deduce, using the monotonicity of $f$ in $y$ that
$$
\ind_{U_s>0}\left( f(s,Y_s,Z_s) - f'\left(s,Y'_s,Z'_s\right) \right) - \mu U^+_s \leq \ind_{U_s > 0} \left( f\left(s,Y'_s,Z_s\right) - f'\left(s,Y'_s,Z'_s\right) \right).
$$
But $ f\left(s,Y'_s,Z'_s\right)  - f'\left(s,Y'_s,Z'_s\right)$ is nonpositive so that
$$
\ind_{U_s>0}\left( f(s,Y_s,Z_s) - f'\left(s,Y'_s,Z'_s\right) \right) - \mu U^+_s \leq \ind_{U_s > 0} \left( f\left(s,Y'_s,Z_s\right) - f\left(s,Y'_s,Z'_s\right) \right).
$$
Finally we set 
$$
\beta_s = \frac{\left(f\left(s,Y'_s,Z_s\right) - f\left(s,Y'_s,Z'_s\right) \right) V_s}{ |V_s|^2 } 
$$
which is a process bounded by $\lambda$.

Coming back to \eqref{tan}, we obtain the following inequality
$$
e^{\mu(t\wedge\tau_n)} U^+_{t\wedge\tau_n} \leq e^{\mu\tau_n} U^+_{\tau_n}    + \int_{t\wedge\tau_n}^{\tau_n} e^{\mu s} \ind_{U_s >0} \beta_s\cdot V_s\, ds - \int_{t\wedge\tau_n}^{\tau_n} e^{\mu s} \ind_{U_s>0} V_s\cdot dB_s 
$$

By Girsanov's theorem, we deduce that 
$$
\e^*\left[ e^{\mu(t\wedge\tau_n)} U^+_{t\wedge\tau_n}\right] \leq \e^*\left[ e^{\mu\tau_n} U^+_{\tau_n} \right],
$$
where $\p^*$ is the probability measure on $(\Omega, \mathcal{F}_T)$ whose density with respect to $\p$ is 
$$
D_T = \exp\left\{\int_0^T \beta_s\cdot dB_s - \frac{1}{2}\int_0^T |\beta_s|^2 ds\right\};
$$
it is worth noting that, since $\beta$ is a bounded process, $D_T$ has moments of all order.

Since we know that $U^+$ belongs to $\m S$, we can easily send $n$ to infinity to get
$$
\e^*\left[ e^{\mu t} U^+_t \right] \leq 0.
$$
Thus $U_t\leq 0$ $\p^*$--a.s. and since $\p^*$ is equivalent to $\p$ on $(\Omega,\mathcal{F}_T)$, $Y_t \leq Y'_t$ $\p$--a.s..
\end{proof}

\subsection{One extension}

In this paragraph, we explain how we can extend our results to a more general setting allowing a superlinear growth of the generator in the variable $y$ as in the work \cite{LsM98}.

Let $h : \rset_+ \fl \rset_+$ be a nondecreasing convex $\m C^1$ function with $h(0)>0$ such that 
$$
\int_0^{+\infty} \frac{du}{h(u)} = +\infty. 
$$

We denote \eqref{h1'} the assumption: there exists $\gamma>0$  such that $\p$--a.s.
\begin{equation}
\label{h1'}
\tag{H1'}
\begin{split}
& \forall t\in[0,T],\qquad (y,z)\longmapsto f(t,y,z) \mbox{ is continuous}, \\
& \forall (t,y,z)\in[0,T]\times\rset\times\rset^d,\qquad \left| f(t,y,z) \right|\leq h(|y|) + \frac{\gamma}{2} |z|^2, \\
& \sup_{y>0} e^{-\gamma y} h(y) < +\infty.
\end{split}
\end{equation}

Let us point out that the previous setting, namely the linear growth condition, corresponds to  $h(y) = \alpha + \beta y$ but we can also have a superlinear growth in $y$; for instance, we can take $h(y) = \alpha\,(y+e) \ln (y+e)$.

Before giving our integrability condition for the terminal value $\xi$, let us explain what is the first modification we have to do. We consider only the case where $h$ is not constant. 

According to the third point of \eqref{h1'}, let us denote by $ c = \sup_{p\in(0,1)} \gamma p h\left(-\frac{\ln p}{\gamma}\right)$ and let us define 
$$
p_0 = \inf\left\{ p\geq 1 : \gamma p h\left(\frac{\ln p}{\gamma}\right) \geq c \right\}.
$$
We define finally 
$$
H(p) = \gamma p h\left(\frac{\ln p}{\gamma}\right) \ind_{p \geq p_0} + c \ind_{p<p_0}. 
$$
Then $H$ is convex and we have the following result.

\begin{lemme}
\label{ode}
Let $z\in\rset$. The differential equation
$$
\phi_t = e^{\gamma z} + \int_t^T H(\phi_s)\, ds,\qquad 0\leq t\leq T,  
$$
has a unique continuous solution $\{ \phi_t(z) \}_{0\leq t\leq T}$ which is decreasing. Moreover, for each $t\in[0,T]$, the map $z\longmapsto \phi_t(z)$ is increasing and continuous. 
\end{lemme}

\begin{proof}
$\phi_t$ is solution  if and only if $u_t = \ln \phi_t /\gamma $ is a solution of the differential equation
$$
u_t' = - \theta(u_t), \quad 0\leq t\leq T, \quad u_T = z \geq 0,
$$
where $\theta(x) = h(x) \ind_{ x \geq \frac{\ln p_0}{\gamma} } + \frac{c}{\gamma} e^{-\gamma x} \ind_{x < \frac{\ln p_0}{\gamma} }$. 
Let us consider the function $\Theta$ defined by
$$
\Theta(x) = \int_{-\infty}^x \frac{1}{\theta(u)}\,du,\quad x\in\rset.
$$
Since $\theta$ is positive, $\Theta$ is an increasing bijection from $\rset$ onto $(0,\infty)$ of class $\m C^1$. It's plain to check that the unique solution to the previous differential equation is $\Theta^{-1}(T-t + \Theta(z) )$ since for any solution we have $\Theta(u_t)' = -1$. Thus 
$$
\phi_t = e^{\gamma \Theta^{-1}(T-t + \Theta(z) )}
$$
and the proof of the lemma is complete.
\end{proof}

We are now in position to give our second assumption.
\begin{equation}
\label{h2'}
\tag{H2'}
\phi_0(|\xi|)\mbox{ is integrable.}
\end{equation}

Exactly as in the linear case, we can prove the following existence result.

\begin{thm}
Let assumptions \eqref{h1'} and \eqref{h2'} hold. Then the BSDE~\eqref{maineq} has at least a solution $(Y,Z)$ such that :
$$
-\frac{1}{\gamma} \ln \e\left(\phi_t(\xi) \:|\: \m F_t \right)\leq Y_t \leq \frac{1}{\gamma} \ln \e\left(\phi_t(\xi)\:|\: \m F_t \right) .
$$
\end{thm}

\section{BSDEs in $\lp^1$}
In this section, we use the method developed before to construct solutions to BSDEs when the data are only integrable. BSDEs with integrable data have been studied in \cite{BDHPS03} and we show that, in the one dimensional case, we can extend the result quoted before. 

Let us recall the framework of \cite{BDHPS03}: assumption (A) holds true for the  random function $f$ if there exist constants $\mu\in\rset$, $\lambda\geq 0$, $\delta\geq 0$ and $\alpha\in(0,1)$ such that
\begin{eqnarray*}
&&(y-y')\cdot \left( f(t,y,z)-f(t,y',z) \right) \leq \mu |y-y'|^2; \\
&&\left| f(t,y,z) - f(t,y,z') \right| \leq \lambda |z-z'|; \\
&& y\fl f(t,y,z) \mbox{ is continuous}; \\
&& \mbox{for each }r>0, \psi_r(t) : = \sup_{|y|\leq r} | f(t,y,0) - f(t,0,0) | \in \lp^1((0,T)\times\Omega); \hspace*{20mm}\\
&& |f(t,y,z)-f(t,y,0)| \leq \delta \left(g_t+|y| + |z| \right)^\alpha,
\end{eqnarray*}
where moreover the progressively measurable processes $(|f(t,0,0)|)$ and $(g_t)$ and the terminal condition $\xi$ satisfy
$$
\e\left[|\xi| + \int_0^T \left( |f(s,0,0)| + g_s \right) ds \right] < +\infty.
$$

Let us recall the following result. 

\begin{lemme}[\cite{BDHPS03}]
Under the assumption {\rm (A)}, BSDE~\eqref{maineq} has a unique solution $(Y,Z)$ such that $Y$ is of class {\rm (D)} and $Z\in \mathrm{M}^\beta$ for some $\beta > \alpha$. Moreover $(Y,Z)\in\mathcal{S}^\beta\times\mathrm{M}^\beta$ for each $\beta \in (0,1)$.
\end{lemme}
We should point out that the result holds true in all dimension not only in the real case.

The last assumption on the generator $f$ does not seem to be very natural; it would be better to have a condition of the type 
$$
|f(t,y,z)| \leq c \left( 1 + |y| + |z|^\alpha\right)
$$
and the remaining of this section is devoted to the construction of a solution under this assumption. But before, we state a comparison result for BSDEs under assumption (A).

\begin{lemme}
\label{compl1}
Let {\rm (A)} holds. Then, if $\xi\leq \xi'$ belong to $\lp^1$, $Y^f_t(\xi)\leq Y^f_t(\xi')$.
\end{lemme}

\begin{proof}
Let us fix $n\in\nset^*$ and denote $\tau_n$ the stopping time
$$
\tau_n = \inf\left\{ t\in[0,T] : \int_0^t \left(|Z_s|^2 + \left|Z'_s\right|^2\right) ds \geq n \right\} \wedge T.
$$

Tanaka's formula leads to the equation, setting $U_t = Y_t-Y'_t$, $V_t=Z_t-Z'_t$, 
\begin{equation*}
\begin{split}
e^{\mu(t\wedge\tau_n)} U^+_{t\wedge\tau_n} \leq & e^{\mu\tau_n} U^+_{\tau_n}   - \int_{t\wedge\tau_n}^{\tau_n} e^{\mu s} \ind_{U_s>0} V_s\cdot dB_s \\
&  + \int_{t\wedge\tau_n}^{\tau_n} e^{\mu s} \left\{ \ind_{U_s >0}\left( f(s,Y_s,Z_s) - f\left(s,Y'_s,Z'_s\right) \right) - \mu U^+_s \right\} ds.
\end{split}
\end{equation*}

We deduce from the previous inequality, using the monotonicity of $f$ in $y$ that
$$
\ind_{U_s>0}\left( f(s,Y_s,Z_s) - f\left(s,Y'_s,Z'_s\right) \right) - \mu U^+_s \leq \ind_{U_s > 0} \left( f\left(s,Y'_s,Z_s\right) - f\left(s,Y'_s,Z'_s\right) \right)
$$
and taking into account the last condition on $f$, the right hand side of the previous inequality is bounded from above by
$$
X_s:=2\delta\left(g_s + \left|Y'_s\right| +\left|Z_s\right| +\left|Z'_s\right|\right)^\alpha .
$$

It follows that
$$
e^{\mu(t\wedge\tau_n)} U^+_{t\wedge\tau_n} \leq  e^{\mu\tau_n} U^+_{\tau_n} + \int_0^T e^{\mu s} X_s\,ds - \int_{t\wedge\tau_n}^{\tau_n} e^{\mu s} \ind_{U_s >0} V_s\cdot dB_s
$$
and thus that
$$
e^{\mu(t\wedge\tau_n)} U^+_{t\wedge\tau_n} \leq  \e\left( e^{\mu\tau_n} U^+_{\tau_n} + \int_0^T e^{\mu s} X_s\,ds \: \Big|\: \mathcal{F}_t \right).
$$
Since $Y$ and $Y'$ belongs to the class (D), we can send $n$ to $\infty$ in the previous inequality (see \cite{BDHPS03} for details) to get 
$$
e^{\mu t} U^+_t \leq  \e\left( \int_0^T e^{\mu s} X_s\,ds \: \Big|\: \mathcal{F}_t \right).
$$
As a byproduct, we deduce that $U^+$ belongs to $\mathcal{S}^p$ as soon as $\alpha p < 1$. Thus we can choose $p>1$ such that $\alpha p <1$.

Since $U^+$ belongs to $\m S$, we can apply Proposition~\ref{wcomp} to conclude the proof.
\end{proof}

From now on we assume that $f$ is continuous and satisfies, for some constants $c\geq 0$ and $\alpha\in(0,1)$,
\begin{equation}
\label{h4} 
\tag{H4}
\left| f(t,y,z) \right| \leq c \left( 1 + |y| + |z|^\alpha \right).
\end{equation}
Moreover we will suppose that $\xi\in\lp^1$. This last assumption is denoted by (H5).

\begin{thm}
Let \eqref{h4} and {\rm (H5)} hold. Then the BSDE~\eqref{maineq} has a solution $(Y,Z)$ such that $Y$ belongs to the class {\rm (D)}. Moreover $(Y,Z)$ belongs to $\m S^\beta\times\mathrm{M}^\beta$ for all $0<\beta<1$.
\end{thm}

\begin{proof}
To prove this result, we use the same approach as in the case of quadratic generators. So let us fix $n\in\nset^*$ and $p\in\nset^*$ and set $\xi^{n,p} = \xi^+\wedge n - \xi^-\wedge p$. Since $f$ is assumed to be a continuous map, we can consider, according to \cite{LsM97}, $\left(Y^f(\xi^{n,p}),Z^f(\xi^{n,p})\right)$ as the minimal solution to the BSDE
$$
Y_t = \xi^{n,p} + \int_t^T f(s,Y_s,Z_s)\,ds - \int_t^T Z_s\cdot dB_s,\quad 0\leq t \leq T.
$$
Since we are dealing with minimal solutions, we have 
$$
Y^f_t(\xi^{n,p+1}) \leq Y^f_t(\xi^{n,p}) \leq Y^f_t(\xi^{n+1,p}),
$$ 
and we set $Y_t = \inf_{p\geq 1} \sup_{n\geq 1} Y^f_t(\xi^{n,p})$.

In order to apply the method described before, we have to find an upper bound independent of $(n,p)$ for $\left| Y^f(\xi^{n,p}) \right|$. For this let us observe that
$$
|f(t,y,z) | \leq g(y,z):=2c \left(1+|y| + |z|^\alpha \wedge |z| \right).
$$

This function $g$ is globally Lipschitz continuous so that we have from the classical comparison theorem
$$
Y^f(\xi^{n,p}) \leq Y^g(\xi^{n,p})
$$
and since we have the same inequality for $-f(t,-y,-z)$,
$$
-Y^f(\xi^{n,p}) \leq Y^g(-\xi^{n,p} ).
$$

But the function $g$ satisfies also the assumption (A) and thus form the comparison theorem in the integrable framework -- Lemma~\ref{compl1} --, we deduce that
$$
\left| Y^f_t(\xi^{n,p} ) \right| \leq Y^g_t(|\xi|).
$$

Let us define $Y_t = \inf_{p\geq 1} \sup_{n\geq 1} Y^f_t(\xi^{n,p})$ and for each $k\geq 1$,
$$
\tau_k = \inf\left\{ t\in[0,T] : Y^g_t(|\xi|) \geq k\right\} \wedge T.
$$
Exactly as in the proof of Theorem~\ref{mainthm}, we construct a process $Z$ such that $(Y,Z)$ solves the BSDE~\eqref{maineq}.

To conclude the proof, let us observe that since $|Y_t| \leq Y^g_t(|\xi|)$, $Y$ belongs to the class (D) and to $\m S^\beta$ for each $\beta\in(0,1)$. It follows from \cite[Lemma 3.1]{BDHPS03} that $Z$ belongs to $\mathrm{M}^\beta$ for $\beta\in(0,1)$.
\end{proof}

%
%

\begin{thebibliography}{BDH{\etalchar{+}}03}

\bibitem[BDH{\etalchar{+}}03]{BDHPS03}
Ph. Briand, B.~Delyon, Y.~Hu, \'E. Pardoux, and L.~Stoica, \emph{L$^p$
  solutions of backward stochastic differential equations}, Stochastic Process.
  Appl. \textbf{108} (2003), 109--129.

\bibitem[Bis78]{Bis78}
J.-M. Bismut, \emph{Contr\^ole des syst\`emes lin\'eaires quadratiques:
  applications de l'int\'egrale stochastique}, S\'eminaire de Probabilit\'es,
  XII (Univ. Strasbourg, Strasbourg, 1976/1977), Lecture Notes in Math., vol.
  649, Springer, Berlin, 1978, pp.~180--264.

\bibitem[EKR00]{ER00}
N.~El~Karoui and R.~Rouge, \emph{Pricing via utility maximization and entropy},
  Math. Finance \textbf{10} (2000), no.~2, 259--276, INFORMS Applied
  Probability Conference (Ulm, 1999).

\bibitem[HIM05]{HIM05}
Y.~Hu, P.~Imkeller, and M.~M{\"u}ller, \emph{Utility maximization in incomplete
  markets}, Ann. Appl. Probab. (2005), To appear.

\bibitem[Kob97]{Kob97}
M.~Kobylanski, \emph{R\'esultats d'existence et d'unicit\'e pour des
  \'equations diff\'erentielles stochastiques r\'etrogrades avec des
  g\'en\'erateurs \`a\ croissance quadratique}, C. R. Acad. Sci. Paris S\'er. I
  Math. \textbf{324} (1997), no.~1, 81--86.

\bibitem[Kob00]{Kob00}
\bysame, \emph{Backward stochastic differential equations and partial
  differential equations with quadratic growth}, Ann. Probab. \textbf{28}
  (2000), no.~2, 558--602.

\bibitem[LSM97]{LsM97}
J.-P. Lepeltier and J.~San~Martin, \emph{Backward stochastic differential
  equations with continuous coefficients}, Statist. Probab. Lett. \textbf{32}
  (1997), no.~4, 425--430.

\bibitem[LSM98]{LsM98}
\bysame, \emph{Existence for {B}{S}{D}{E} with superlinear-quadratic
  coefficient}, Stochastics Stochastics Rep. \textbf{63} (1998), no.~3-4,
  227--240.

\bibitem[PP90]{PP90}
E.~Pardoux and S.~Peng, \emph{Adapted solution of a backward stochastic
  differential equation}, Systems Control Lett. \textbf{14} (1990), no.~1,
  55--61.

\end{thebibliography}
%
\newcommand{\etalchar}[1]{$^{#1}$}
\providecommand{\bysame}{\leavevmode\hbox to3em{\hrulefill}\thinspace}
\providecommand{\MR}{\relax\ifhmode\unskip\space\fi MR }
\providecommand{\MRhref}[2]{%
  \href{http://www.ams.org/mathscinet-getitem?mr=#1}{#2}
}
\providecommand{\href}[2]{#2}

\end{document}